\documentclass[11pt]{article}

\usepackage{graphicx}
\usepackage{amsmath}
\usepackage{amssymb}
\usepackage{epstopdf}
\usepackage[all]{xy}
\usepackage{tensor}

\renewcommand \baselinestretch 2

\newcommand {\m}[1]{}
\newcommand{\bbrc}{\left(\begin{array}{c}}
\newcommand{\ebrc}{\end{array}\right)}
\newcommand{\omg}{ {\omega,g} }

\newtheorem{definition}{Definition}
\newtheorem{prop}{Proposition}
\newtheorem{lem}{Lemma}
\newtheorem{rem}{Remark}

\newenvironment{pf}{{\it proof: }}{\hfill $\spadesuit$}

\begin{document}

\newcommand{\gan}{\rangle}
\newcommand{\lan}{\langle}

\title{On two isomorphic Lie algebroids for Feedback Linearization}

\author{M\"ullhaupt Philippe  \\ \\
D\'epartement de g\'enie m\'ecanique, EPFL, CH-1015 Lausanne}

\maketitle

\begin{abstract}
Two Lie algebroids are presented that are linked to the construction of the linearizing output of an affine in the input nonlinear system. 
 The algorithmic construction of the linearizing output proceeds inductively, and each stage has two structures, namely a codimension one foliation defined through an integrable 1-form $\omega$ , and a transversal vectorfield $g$  to the foliation. Each integral manifold of the vectorfield $g$ defines an equivalence class of points. Due to transversality, a leaf of the foliation is chosen to represent these equivalence classes. A Lie groupoid is defined with its base given as the particular chosen leaf and with the product induced by the pseudogroup of diffeomorphisms that preserve equivalence classes generated by the integral manifolds of g. Two Lie algebroids associated with this groupoid are then defined. The theory is illustrated with an example using polynomial automorphisms as particular cases of diffeomorphisms and shows the relation with the Jacobian conjecture.
\end{abstract}



\section{Introduction}

Affine in the input nonlinear systems (\cite{Isidori}, \cite{Nij}) are considered with a single control $u$  and with state  $x \in \mathbb{R}^n$ defined by
\begin{eqnarray}
\dot x &=& f(x) + g(x) u \nonumber 
\end{eqnarray}
This system is feedback linearizable to a linear system $\dot z = A z + B v$ through diffeomorphism $z = \Phi(x)$ and change of coordinates $v = \alpha(x) + \beta(x) u$ under the condition of accessibility, i.e. $\mbox{rank}(g, \mbox{ad}_f g, \ldots, \mbox{ad}_f^n)$ and involutivity of the distribution ${\mathcal C} = \mbox{span}\{(g, \mbox{ad}_f g, \ldots, \mbox{ad}_f^{n-2} g\}$ (\cite{Isidori}, \cite{Nij}). A classical way of computationally solving this problem is to use the flow-box theorem \cite{Tall} which amounts to inductively straighten out the vectorfields. A similar method is used in the proof of the Frobenius theorem in \cite{Chevalley} Theorem 9 on pp. 89-92, and in \cite{BishopCrittenden}, Theorem 7 on p. 24. Another approach is to integrate the integrable 1-form in the null-space of the distribution $\mathcal C$ and relates to the dual approach of \cite{GardnerShadwick}, \cite{Schlacher}, \cite{Basile}. Equivalence in the classical setting between the two approaches can be found on p. 71 of \cite{BishopCrittenden}.

An inductive process using a somewhat intermediate approach between the two appeared in \cite{Mullhaupt} where an anti-symmetrical product was defined.

The point of the following developments is to throw light on the meaning of the anti-symmetrical product defined in \cite{Mullhaupt} by proving that it is a Lie algebroid. This is achieved through a tedious albeit direct proof of the Jacobi identity and the definition of a suitable anchor map. 
In \cite{Mullhaupt}, this Lie algebroid was related to a Lie groupoid without mentioning this formalism. 

In \cite{Willson} another anchor map was defined without explicitly mentioning the Lie algebroid formalism. Clarification of the relations between the two algebroids (by providing an isomorphism of algeboroids) and between the algebroids and the groupoid will be given.

	An interesting application of the theory is provided when the diffeomorphism of the definition of feeback linearization is replaced by a polynomial automorphism (see \cite{vanEssen} for a detailed coverage of this topic in relation with the Jacobian conjecture). The intermediate 1-forms appearing in the definition of the algebroid when suitably defined leads to an algorithm for finding the polynomial inverse map of the polynomial automorphism $z = \Phi(x)$. If all the 1-forms appearing throughout the intermediate steps (where the anchor map is used) could be shown to have have constant determinant, this would lead to the proof of the Jacobian Conjecture. 

Section \ref{LieAlgebroidGroupoid} introduced the definition of a Lie groupoid of the literature, fixes notations, and gives explicitly the axioms for the class of Lie groupoids that will be used with feedback linearization. We also recall the definition of a Lie algebroid and define the two aforementioned Lie algebroids. The proof of the Jacobian identity is then given for the first algebroid together with an inductive construction of the linearizing output using Algebroid I and Algebroid II. Section \ref{AppPolAuto} applies the theory to the case of polynomial automorphisms and relates both algorithms to the Jacboian Conjecture. Complete proofs omitted due to the page limit can be found in \cite{MullhauptExtendedVersion}.

\section{Lie Groupoid and Lie Algebroid}

\label{LieAlgebroidGroupoid}

\subsection{Lie Groupoid}

A lie groupoid \cite{Kiril2}, \cite{Moerdijk} consists of six elements subject to five axioms. 

\begin{definition}{\sc Lie Groupoid. }
\label{defGroupoid}
A Lie groupoid \cite{Kiril2}, \cite{Moerdijk} consists of the six elements: 
\begin{enumerate}
\item[I.] A set $\Omega$ called the groupoid (set of arrows)
\item[II.] a set $\mathcal O$ called the base (set of objects)
\item[III.] a source map $\sigma$, from $\Omega$ to $\mathcal O$
\item[IV.] a target map $\tau$, from $\Omega$ to $\mathcal O$
\item[V.] an object inclusion map $\iota$, from $\mathcal O$ to $\Omega$
\item[VI.] a partial multiplication map $(\Phi_1, \Phi_2) \rightarrow  \Phi_1 \bot \Phi_2$, from $\Omega * \Omega$ to $\Omega$, where 
$$\Omega * \Omega = \{Ê(\Phi_1, \Phi_2) \in \Omega \times \Omega \; | \; \sigma(\Phi_1) = \tau(\Phi_2) \}$$
\end{enumerate}
The target map and the source map are surjective submersions. The inclusion map is smooth. The partial multiplication $\bot$ is smooth. Additionally,  
these six elements are subject to the axioms:
\begin{enumerate}
\item[(i)] $\sigma(\Phi_1 \bot \Phi_2) = \sigma(\Phi_2)$ and $\tau(\Phi_1 \bot \Phi_2) = \tau(\Phi_1)$ for all $(\Phi_1, \Phi_2) \in \Omega * \Omega$;
\item[(ii)] $\Phi_1 \bot (\Phi_2 \bot \Phi_3) = (\Phi_1 \bot \Phi_2)Ê\bot \Phi_3$ for all $\Phi_1, \Phi_2, \Phi_3 \in \Omega$ such that $\sigma(\Phi_1) = \tau(\Phi_2)$ and $\sigma(\Phi_2) = \tau(\Phi_3)$;
\item[(iii)] $\sigma(\iota(\bar O)) = \tau(\iota(\bar O)) = \bar O$ for all $\bar O \in \mathcal O$;
\item[(iv)] $\Phi_2 \bot \iota(\sigma(\Phi_2)) = \Phi_2$ and $\iota(\tau(\Phi_2)) \bot \Phi_2 = \Phi_2$ for all $\Phi_2 \in \Omega$;
\item[(v)] each $\Phi_2 \in \Omega$ has an inverse $\Phi_2^{-1}$ such that 
\begin{itemize}
\item[]
$\sigma(\Phi_2^{-1}) = \tau(\Phi_2)$,
$\tau(\Phi_2^{-1}) = \sigma(\Phi_2)$
\item[] 
 $\Phi_2^{-1} \bot \Phi_2 = \iota(\sigma(\Phi_2))$, $\Phi_2 \bot \Phi_2^{-1} = \iota(\tau(\Phi_2))$. 
 \end{itemize}
\end{enumerate}
The element $\iota(\bar O)  \in \Omega$ corresponding to ${\bar O} \in \mathcal O$ may be called the unity or identity corresponding
to $\bar O$. 
\end{definition}

\subsection{The Lie Groupoid for Feedback Linearization}

\label{Groupoid}
A vectorfield $g$ is given together with a noncancelling integrable 1-form $\omega$, that is, $\omega g \neq 0$ 
for all $x \in \mathbb{R}^n$ and $d \omega \wedge \omega = 0$, where $d$ stands for the exterior derivative. 
This means that $\omega$ admits locally integral manifolds constituting a codimension $1$ foliation (see for example \cite{Lawson}).

\begin{definition}
An integral manifold of $\omega$ passing through a point $A$ of the surrounding manifold will be written as
$\mathcal{O}_A$.
\end{definition}

Because the distribution defined by the vectorfield $g$ is trivially involutive and nonvanishing, it admits integral manifolds:

\begin{definition}
The integral manifold of the vectorfield $g$ passing through a point $A$ of the surrounding manifold
is designated by $\mathcal{G}_A$.
\end{definition}

Lemma \ref{integLem} shows that the set of all diffeomorphisms preserve the foliation defined by $\omega$, since $\omega$ is assumed integrable. The groupoid under study will be a subset of these
diffemorphisms that preserve equivalence classes defined by integral manifolds $\mathcal G$ of $g$.   

\begin{definition}
\label{equivClass}{\sc Equivalence classes along integral manifolds of $g$}
Two points $A_1$ and $A_2$ belong to the same equivalence class 
whenever
$$A_1 \in {\mathcal G}_{A_2},$$ or, what means the same thing, whenever
$$A_2 \in {\mathcal G}_{A_1}.$$  
\end{definition}

\begin{definition}{\sc Elements $\Omega_I$. }
\label{OmegElemsEquiv}
\label{defOmegI}
Elements of $\Omega_I$ are diffeomorphisms $\Phi_{A,B}$ such that:
\begin{itemize}
\item they map the point $A$ to the point $B$, i.e. $\Phi_{A,B}(A) = B$;
\item they preserve integral manifolds of $g$:
$$\forall C \in {\mathcal G}_A \cap \mathfrak{D}(\Phi_{A,B}) \Rightarrow \Phi_{A,B}(C) \in \Phi_{A,B}({\mathcal G}_A) \cap \mathfrak{R}(\Phi_{A,B}).$$
\end{itemize}
\end{definition}

\begin{definition}{\sc Elements $\Omega_{II}$.Ê}
\label{defOmegII}
Let $\psi_j: \mathbb{R}^{n-1} \rightarrow \mathbb{R}$, $j=A,B$ be two functions satisfying 
both $\psi_j(j) = 0$, $j=A,B$ and $d \psi_j = \mu_j \omega$, $j=A,B$ with two functions
 $\mu_j : \mathbb{R}^n \rightarrow \mathbb{R}$.
Choosing $n-1$ functions $\phi_{A,i}$, $i=1, \ldots, n-1$ such that (i) $\phi_{A,i}(A) = 0$, $i=1, \ldots, n-1$ and (ii) the 1-forms $d \phi_{A,i}$, $i=1, \ldots, n-1$ together with $\omega$, evaluated at $A$, constitute a basis
of $T_A^* \mathbb{R}^n$ and (iii)
$d \phi_{A,i} g = 0$, $i=1, \ldots, n-1$.
Similarly, choose another set of functions $\phi_{B,i}$, $i=1, \ldots, n-1$,
so that (i) $\phi_{B,i}(A) = 0$, $i=1, \ldots, n-1$ and (ii) $d \phi_{B,i}$, $i=1, \ldots, n-1$ together with $\omega$, evaluated at $B$, constitute a basis
of $T_B^* \mathbb{R}^n$ and (iii) $d \phi_{B,i} g = 0$, $i=1, \ldots, n-1$. Then $\Omega_{II}$ is the set of all diffeomorphisms $\Phi_{A,B} : \mathbb{R}^n \rightarrow
\mathbb{R}^n$ that can be expressed as 
\begin{equation}
\label{phiCoords}
\Phi_{A,B} := \Phi_B^{-1}  \circ \Phi_A
\end{equation}
with 
\begin{equation}
\label{phidef}
\Phi_A := \left(\begin{array}{c} \phi_{A,1} \cr \phi_{A,2} \cr \vdots \cr \phi_{A,n-1} \cr \psi_A \end{array} \right) \qquad 
\Phi_B := \left(\begin{array}{c} \phi_{B,1} \cr \phi_{B,2} \cr \vdots \cr \phi_{B,n-1} \cr \psi_B \end{array} \right) 
\end{equation}
\end{definition}

\begin{lem}
The set $\Omega_{II}$ is a subset of $\Omega_I$.
\end{lem}

\begin{pf}
Because the corresponding constituting 1-forms $d \psi_A$, $d \phi_{A,1}$, $d \phi_{A,2}$, $\ldots$, $d \phi_{A,n-1}$ (resp. $d \psi_B$, $d \phi_{B,1}$, $d \phi_{B,2}$, $\ldots$, $d \phi_{B,n-1}$) form a basis of $T_A^{*} \mathbb{R}^n$ (resp. $T_B^{*} \mathbb{R}^n$), when evaluated at $A$ (resp. $B$), the maps $\Phi_A$ and $\Phi_B$ in (\ref{phidef}) are local diffeomorphisms, so that the reciprocal map $\Phi_B^{-1}$ exists showing that 
(\ref{phiCoords}) is a well defined diffeomorphism. Additionally, $\mathfrak{D}(\Phi_{A,B}) = \mathfrak{R}(\Phi_{A,B}) = \mathbb{R}^n$.
Let $x$ designate the coordinates of the surrounding manifold $\mathbb{R}^n$. 
Define $z$-coordinates as $z_n := \phi_1(x)$, $z_2 := \phi_{A,2}(x)$, $\ldots$, $z_{n-1} := \phi_{A,n-1}(x)$,
$z_n := \psi_A(x)$. Then set ${\mathcal O}_A := \{Êx | \psi_A(x) = 0\}$ so that ${\mathcal O}_A$ is both a local integral manifold of $\omega$ and a set that contains $A$. In the $z$ coordinates, its expression is ${\mathcal O}_A = \{Êz |Êz_n = 0\}$.
Similarly, define $z_1' := \phi_{B,1}(x)$, $z_2' := \phi_{B,2}(x)$, $\ldots$, $z_{n-1}' := \phi_{B,n-1}(x)$,
$z_n' := \psi_B(x)$ so that setting ${\mathcal O}_B := \{ x | \psi_B(x) = 0 \}$ defines both a local integral manifold 
of $\omega$ and a set containing $B$. Expressed in the $z'$ coordinates, ${\mathcal O}_B = \{Êz' |Êz_n' = 0\}$.
Now, the choices (\ref{phidef}) defining (\ref{phiCoords}) show that
the composition operator appearing in (\ref{phiCoords}) forces $z_n = z_n'$ so that 
$\Phi_{A,B}( {\mathcal O}_A) = {\mathcal O}_B$ which confirms that $\Phi_{A,B} \in \Omega_I$ according to
Definition \ref{OmegElemsEquiv}.
\end{pf}

\begin{definition}{\sc Base manifold $\mathcal O$. }
The base manifold $\mathcal O$ is a globally defined integral manifold of $\omega$.
\end{definition}

\begin{definition}{\sc Function $\psi$. }
\label{psiFunction}
We will suppose that $\mathcal O$ is defined by a single function $\psi : \mathbb{R}^n \rightarrow \mathbb{R}$ through \begin{equation}
\label{baseSet}
{\mathcal O} = \{x | \psi(x) = 0\}.
\end{equation}
\end{definition}

\begin{definition}
\label{sourceMap}{\sc Source map $\sigma$. } The source map $\sigma$ maps the domain $\mathfrak{D}( \Phi_{A,B})$ of 
a diffemorphism $\Phi_{A,B} \in \Omega$ to the base manifold $\mathcal O$ by following integral manifolds $\mathcal G$ of $g$, that is,
$$\sigma(A_1) := {\mathcal G}_{A_1} \cap \mathcal O, \quad \forall A_1 \in \mathfrak{D}(\Phi_{A,B}).$$
\end{definition}

\begin{rem}
Notice that Definition \ref{sourceMap} is well defined because we assume $\omega g \neq 0$ globally. The groupoid can be understood as a class of pseudo-group. Pseudo-groups are used when dealing with accessible sets \cite{Stefan} and with Riemannian foliations \cite{Molino}.
\end{rem}

\begin{definition}{\sc Target map $\tau$. }
\label{targetMap} The target map $\tau : \mathbb{R}^n \rightarrow \mathbb{R}$ maps the range $\mathfrak{R}( \Phi_{A,B})$ of 
an element $\Phi_{A,B}$ to the base manifold $\mathcal O$ by following integral manifolds $\mathcal G$ of
$g$, that is,
$$\sigma(B_1) := {\mathcal G}_{B_1} \cap \mathcal O, \quad \forall B_1 \in \mathfrak{R}(\Phi_{A,B}).$$
\end{definition}

\begin{lem}
\label{phiLemma}
Under the hypothesis of the existence of a function $\psi$ according to Definition \ref{psiFunction} and of the existence of a base of 1-forms of $T^* \mathbb{R}^n$, both the source map $\sigma$ (Definition \ref{sourceMap}) and the target map $\tau$ (Definition \ref{targetMap}) are globally defined and can be described using coordinates by choosing 
$n-1$ functions $\gamma_1$, $\gamma_2$, $\ldots$, $\gamma_{n-1}$ such that $d \gamma_i g = 0$, $i=1, \ldots, n-1$
and such that $d \gamma_i$, $i=1,\ldots, n-1$ are independent 1-forms.
\end{lem}

\begin{pf}
Because the corresponding constituting 1-forms $d \psi_A$, $d \phi_{A,1}$, $d \phi_{A,2}$, $\ldots$, $d \phi_{A,n-1}$ (resp. $d \psi_B$, $d \phi_{B,1}$, $d \phi_{B,2}$, $\ldots$, $d \phi_{B,n-1}$) form a basis of $T_A^{*} \mathbb{R}^n$ (resp. $T_B^{*} \mathbb{R}^n$), when evaluated at $A$ (resp. $B$), the maps $\Phi_A$ and $\Phi_B$ in (\ref{phidef}) are local diffeomorphisms, so that the reciprocal map $\Phi_B^{-1}$ exists showing that 
(\ref{phiCoords}) is a well defined diffeomorphism. Additionally, $\mathfrak{D}(\Phi_{A,B}) = \mathfrak{R}(\Phi_{A,B}) = \mathbb{R}^n$.
Let $x$ designate the coordinates of the surrounding manifold $\mathbb{R}^n$. 
Define $z$-coordinates as $z_n := \phi_1(x)$, $z_2 := \phi_{A,2}(x)$, $\ldots$, $z_{n-1} := \phi_{A,n-1}(x)$,
$z_n := \psi_A(x)$. Then set ${\mathcal O}_A := \{Êx | \psi_A(x) = 0\}$ so that ${\mathcal O}_A$ is both a local integral manifold of $\omega$ and a set that contains $A$. In the $z$ coordinates, its expression is ${\mathcal O}_A = \{Êz |Êz_n = 0\}$.
Similarly, define $z_1' := \phi_{B,1}(x)$, $z_2' := \phi_{B,2}(x)$, $\ldots$, $z_{n-1}' := \phi_{B,n-1}(x)$,
$z_n' := \psi_B(x)$ so that setting ${\mathcal O}_B := \{ x | \psi_B(x) = 0 \}$ defines both a local integral manifold 
of $\omega$ and a set containing $B$. Expressed in the $z'$ coordinates, ${\mathcal O}_B = \{Êz' |Êz_n' = 0\}$.
Now, the choices (\ref{phidef}) defining (\ref{phiCoords}) show that
the composition operator appearing in (\ref{phiCoords}) forces $z_n = z_n'$ so that 
$\Phi_{A,B}( {\mathcal O}_A) = {\mathcal O}_B$ which confirms that $\Phi_{A,B} \in \Omega_I$ according to
Definition \ref{OmegElemsEquiv}.
\end{pf}

\begin{definition} {\sc $\Phi$ map.}

Define the $\Phi$ map
as
\begin{equation}
\Phi = \begin{pmatrix} 
\gamma_1(x_1, \ldots, x_n) \\
\vdots \\
\gamma_{n-1}(x_1, \ldots, x_n) \\
\psi(x_1, \ldots, x_n) 
\end{pmatrix}
\end{equation}
so that according to Lemma \ref{phiLemma} both the source map and target map can be defined as $\sigma = \Phi$ and $\tau = \Phi$.
\end{definition}

\begin{definition}{\sc Inclusion map $\iota$.} The inclusion map $\iota({\bar B})$ associates a diffeomorphism $\Phi_{B,B}: \mathbb{R}^n \rightarrow \mathbb{R}^n$ to the the point $\bar B \in {\mathcal O}$, with $B$ being the inclusion of $\bar B$ in the surrounding manifold $\mathbb{R}^n$, such that  
$\Phi_{B,B}$ is an identity on a  
local submanifold ${\mathcal O}_{\bar B}$ of $\omega$ (of same dimension) that contains $\bar B$.
\end{definition}

\begin{definition}{\sc Product $\bot$}
Given two elements $\Phi_{A_1, B_1}$ and $\Phi_{B_2, C_2}$ of $\Omega_I$ for which $B_1 \in \mathcal{G}_{B_2}$, define their product as
\begin{equation}\label{botProductI}\Phi_{A_1, B_1} \bot  \Phi_{B_2, C_2} := \Phi_{B_2, C_2} \circ \Phi_{A_1, B_1}.
\end{equation} 
\end{definition}

\begin{prop}
Axioms (i) to (v) of a Lie groupoid appearing in Definition \ref{defGroupoid} are satisfied for elements of $\Omega_I$
given in Definition \ref{defOmegI} and for the product (\ref{botProductI}).
\end{prop}

\begin{pf}
Axiom $(i)$ is satisfied by definition of $\Phi_{A_1, C_1}$ because it shares the same $\alpha$ map, i.e.
$\alpha_{A_1}$ for $\Phi_{A_1, C_1}$ is the same as $\alpha_{A_1}$ for $\Phi_{A_1, B_1}$.
Axiom (ii) is trivially satisfied because of the associativity of compositions of maps.
The object inclusion map $\iota$ is the identity map
$$\iota : B \rightarrow \mathbb{R}^n \cap \{Êx \; | \; \psi(x) = 0 \}$$
so that Axiom (iii), which is $\alpha(\iota(\bar O)) = \beta(\iota(\bar O))$, is also satisfied.
However, Axiom (iv) is slightly more involved. Let us suppose that
$\xi = \Phi_{A_1, B_1}$ so that $\xi$ maps $O_{A_1}$ to $O_{B_1}$.  Then $\sigma(\Phi)$ is the map between
$O_{A_1}$ to $\mathcal O$ that assigns to every point of $A \in O_{A_1}$ the point ${\mathcal G}_{A} \cap {\mathcal O}$ in $\mathcal O$.
Therefore, if one mutiplies by  $\Phi$, that is  $\Phi \bot \iota(\sigma(\Phi))$, then one gets back $\Phi$ because
of the correspondence along the integral manifolds of $g$ between the image of $O_{A_1}$ as an open set in $\mathcal O$  and $O_{A_1}$ itself. 
\end{pf}

\subsection{Lie Algebroid}

\begin{definition}{\sc Lie Algebroid. }
\label{AlgebroidDefinition}
Let $\mathcal O$ be a manifold. A Lie algebroid on $\mathcal O$ is a vector bundle $(A, \pi, {\mathcal O})$ together with a vector bundle map
$\pi : A \rightarrow T {\mathcal O}$ over $\mathcal O$, called the anchor of $A$, and a bracket on sections $\Gamma A$
of the bundle given as $[.,.] : \Gamma A \times \Gamma A \rightarrow \Gamma A$ which is $\mathbb{R}$-bilinear and alternating
$$[m_1,m_2]  = - [m_2,m_1] \qquad \qquad m_1, m_2  \in \Gamma A$$ and satisfies Jacobi's identity, i.e. $\forall m_1, m_2, m_3  \in \Gamma A$,
$$[m_1,[m_2,m_3]] + [m_2,[m_3,m_1]] + [m_3,[m_1,m_2]] = 0$$ 
The anchor and the bracket satisfy the properties:
\begin{enumerate}
\item[(I)] $\pi([m_1, m_2]) = [\pi(m_1), \pi(m_2)] \qquad \qquad m_1, m_2  \in \Gamma A$\item[(II)] $[m_1, \alpha m_2] = \alpha [m_1, m_2] + \left( L_{\pi(m_1)} \alphaÊ\right) m_2 \qquad m_1, m_2 \in \Gamma A, \alpha \in C({\mathcal O})$. 
\end{enumerate}
where $C({\mathcal O})$ designates functions on $\mathcal O$.
\end{definition}

\subsection{Effect of diffeomorphisms on vectorfields and 1-forms}

Consider an arbitrary diffeomorphism $\Phi: \mathbb{R}^n \rightarrow \mathbb{R}^n$. Using coordinates, $\Phi$
defines a new set of coordinates $z$ using the initial coordinates $x$ as $z := \Phi(x)$.
This has consequences on vectorfields belonging to $T \mathbb{R}^n$ and 1-forms belonging to $T^*\mathbb{R}^n$.

\begin{definition}{\sc Push-forward.}
Let $m \in T \mathbb{R}^n$ be a vectorfield. Define the push-forward of $m$ by the diffeomorphism $\Phi$ by 
\begin{equation}
\label{pushforward}
\Phi_*(m) := \frac{\partial \Phi}{\partial x} m \circ \Phi^{-1}(z)
\end{equation}
\end{definition}

\begin{definition}{\sc Pull-back.}
Let $\omega \in T^*\mathbb{R}^n$ be a 1-form. Using the vector notation that associates to the 1-form
$\sum_{i=1}^n \omega_i(x) dx_i$ the vector $\omega = \begin{pmatrix} \omega_1 & \omega_2 & \ldots & \omega_nÊ\end{pmatrix}$, 
define the pull-back of $\omega$ by $\Phi$ by
$$\Phi^*(\omega) := \omega \left(\frac{\partial \Phi}{\partial x}\right)^{-1} \circ \Phi^{-1}(z)$$
\end{definition}

\begin{lem}
If $m$ is a tangent vector to a curve ${\mathcal C} = \{ x | x = \xi(\alpha), \alpha \in \mathbb{R}\}$ with
$\xi : \mathbb{R} \rightarrow \mathbb{R}^n$ a smooth defining funtion, then $\Phi_*(m)$ is the tangent 
vector of the image $\Phi({\mathcal C}) := \{ z | z = \Phi(\xi(\alpha)), \alpha \in \mathbb{R} \}$ of the curve $\mathcal C$ under the diffeomorphism $\Phi$.
\end{lem}

\begin{lem}
\label{integLem}
If $\omega$ is an integrable 1-form associated with the integral manifold locally defined by a function $\psi:  \mathbb{R}^n \rightarrow \mathbb{R}$ as $\{ x | \psi(x) = 0\}$,
then the pull-back $\Phi^* \omega$ remains an integrable 1-form. Moreover, $\psi \circ \Phi^{-1}$ defines
locally an integral manifold of $\Phi^* \omega$. This manifold is locally described as the set $\{ z | \psi \circ \Phi^{-1}(z) = 0 \}$. 
\end{lem}

\begin{pf}
These two results are classical, see for example \cite{Morita}. 
\end{pf}

\subsection{Lie Algebroid I for Feedback Linearization}

The bracket is defined as
\begin{eqnarray}
\label{equivClass}
\lan \bar m_1, \bar m_2 \gan \simeq \lan m_1, m_2 \gan  \nonumber \\
:= [m_1, m_2] + \frac{\omega m_2}{\omega g} [g, m_1] - \frac{\omega m_1}{\omega g} [g, m_2]
\end{eqnarray}
where $m_1$ (resp. $m_2$) is any representative of the equivalence class of $\bar m_1$ (resp. $\bar m_2$). 
This definition of the  anti-symetrical product appeared in \cite{Mullhaupt} without either the Lie algebroid interpretation or mentioning the equivalence classes on which it operates. The closest definition that the author could find is the Nickerson bracket, i.e.   formula (44) on p. 520 in \cite{Nickerson}. The explicit appearance of the integrable 1-form $\omega$ does however not appear in that formula.

\begin{lem}
The bracket in (\ref{equivClass}) is independent of the equivalence classes $m_1$ and $m_2$ chosen.
\end{lem}

\begin{pf}
\begin{eqnarray}
&&\lan \bar m_1, \bar m_2 \gan \simeq \lan m_1 + \alpha g, m_2 + \beta g \gan  \nonumber \\
&=& [m_1 + \alpha g, m_2 + \beta g] + \frac{\omega (m_2+\beta g)}{\omega g}[g , m_1 + \alpha g] \nonumber \\
&& - \frac{\omega(m_1 + \alpha g)}{\omega g} [g, m_2] \nonumber \\
&=& [m_1,m_2] + \beta [m_1,g] + \alpha [g, m_2] \nonumber \\
&&+ \left(m_1(\beta) - m_2(\alpha) + \alpha g(\beta) - \beta g(\alpha) \right) g \nonumber \\
&& + \frac{\omega m_2}{\omega g} [g, m_1] + \frac{\omega m_2}{\omega g} g(\alpha) g + \beta [g,m_1] + \beta g(\alpha) g \nonumber \\
&& - \frac{\omega m_1}{\omega g} [g, m_2] - \frac{\omega m_1}{\omega g} g(\beta) g - \alpha [g, m_2] - \alpha g(\beta) g \nonumber \\
&=& \lan m_1, m_2 \gan + \left( m_1(\beta) - m_2(\alpha) + \frac{\omega m_2}{\omega g} g(\alpha) - \frac{\omega m_1}{\omega g} g(\beta)  \right)g \nonumber \\
&\simeq& \lan \bar m_1, \bar m_2 \gan 
\end{eqnarray}
\end{pf}

\subsection{Lie Algebroid on $(\mathcal{O}, T\mathbb{R}^2/\mathcal{G})$}

\label{Algebroid1}

The base manifold $\mathcal O$ is an integral manifold of the integrable 1-form $\omega \in T^*{\mathbb{R}^n}$ and the typical fibre bundle is $T{\mathbb{R}^n}_x/\mbox{span }g(x)$, a section of which  is a map 
$m : \mathcal{O} \rightarrow T\mathbb{R}^2/\mathcal{G}$.

\subsubsection{The Anchor}

\begin{definition}
Let $\mathcal{O}$ designate an integral manifold of the integrable 1-form $\omega$.
The following anchor $\mbox{an }_\pi: T\mathbb{R}^n \rightarrow  T{\mathcal O}$ is defined as
\[
\mbox{an }_\pi(m) := \pi_{\omega,g \,*} m
\]
where $\pi_{\omega,g}$ is the  projection operator $\pi_{\omega,g} : \mathbb{R}^n \rightarrow {\mathcal O}$ along integral curves of $\mathcal G$, i.e. $\pi_{\omega, g}(m_1) = \pi_{\omega,g}(m_2)$ whenever $m_1 \in {\mathcal G}_{m_2}$ (i.e. $m_2 \in {\mathcal G}_{m_1}$). It is such that $\pi_{\omega,g \, *}(g) =0$.
\end{definition}

\subsection{Properties I and II of the anchor $\mbox{an }_\pi$}

\begin{lem}
With anchor $\mbox{an }_\pi$, Property I holds:
\[
\lan \bar m_1, \alpha \bar m_2 \gan = \alpha \lan \bar m_1, \bar m_2 \gan + \mbox{an }_\pi(\bar m_1)(\alpha) \bar m_2 \qquad \forall \alpha \in C({\mathcal O})
\]
\end{lem}

\begin{pf}
The function $\alpha \in C({\mathcal O})$ can be expressed with coordinates $z_1, \ldots, z_{n-1}$ that locally defines the embedded submanifold $\mathcal O$. Hence we can also understand $\alpha$ as defined in $\mathbb{R}^n$ by considering $\alpha$ as a function of $z_1, \ldots, z_n$ with $z_n =0$ defining $\mathcal O$. Denote the change of coordinates from $x$ in $\mathbb{R}^n$ to $z$ by $z = \Phi(x)$. This then means that $\Phi_* g = \frac{\partial}{\partial z_n}$ by construction of $\pi_{\omega, g, \; *} = \mbox{Pr }\Phi_* g$ where $\mbox Pr$ meaning the projection by not considering the last coordinate. Since $\alpha$ does not depend on $z_n$ by construction, it holds that $L_g \alpha = 0$, so that
\begin{eqnarray}
&&\lan m_1, \alpha m_2 \gan \nonumber \\
 &=& [m_1, \alpha m_2] + \frac{\omega \alpha m_2}{\omega g} [g, m_1] - \frac{\omega m_1}{\omega g}[g, \alpha m_2] \nonumber \\
&=& \alpha [m_1, m_2] + m_1(\alpha) m_2  \nonumber \\
&&+ \alpha \frac{\omega m_2}{\omega g} [g, m_1] - \alpha \frac{\omega m_1}{\omega g} [g, m_1] - \frac{\omega m_1}{\omega g} g(\alpha) m_2 \nonumber \\
&=& \alpha \lan m_1, m_2 \gan + m_1(\alpha) m_2 \nonumber
\end{eqnarray}
Now since $g(\alpha) = 0$, it follows that $m_1(\alpha) = \pi_{g, \omega *} m_1(\alpha) = \mbox{an }_\pi(m_1)(\alpha)$ proving the required identity.
\end{pf}

\begin{lem}
With anchor $\mbox{an }_\pi$, Property II holds:
\begin{equation}
\label{RHS2}
\mbox{an }_\pi( \lan \bar m_1, \bar m_2 \gan) = [ \mbox{an }_\pi(\bar m_1), \mbox{an }_\pi(\bar m_2)]
\end{equation}
\end{lem}

\begin{pf}
The lemma and its proof are given in \cite{Mullhaupt}, Lemma 1 at the bottom of p. 554.  \end{pf}

\subsection{Lie Algebroid on the bundle $(\mathbb{R}^n, \mathbb{R}^n/{\mathcal G})$}

\label{Algebroid2}

The base manifold $\mathcal{O}$ is an integral manifold of the integrable 1-form $\omega \in {T\mathbb{R}^n}^*$ and the typical fibre bundle is $T\mathbb{R}^n_x/\mbox{span }g(x)$, for which a section is a map 
$m : \mathbb{R}^n \rightarrow T\mathbb{R}^n§/\mathcal{G}$.

\subsubsection{The Anchor} 

\begin{definition}
Then anchor $\mbox{an }_{\omega,g}: T\mathbb{R}^n/{\mathcal G} \rightarrow T\mathbb{R}^n$ is defined for any any 1-form $\omega$ such that $\omega g \neq 0$. For a given section $\bar m \in \Gamma T\mathbb{R}^n/{\mathcal G}$, the anchor is defined as
\begin{equation}
\mbox{an }_{\omega,g}(\bar m) := m - \frac{\omega m}{\omega g} g
\end{equation}
where $m$ is any representative in $\Gamma T\mathbb{R}^n$ of the equivalence class $\bar m \in T\mathbb{R}^n/{\mathcal G}$.
\label{anchorWG}
\end{definition}

\begin{lem}
The elements in Definition \ref{anchorWG} are well defined
\end{lem}

\subsubsection{Properties I and II of the anchor $\mbox{an }_{\omega,g}$}

\begin{lem}
Property I holds:  
\[
\lan \bar m_1, \alpha \bar m_2 \gan = \alpha \lan \bar m_1, \bar m_2 \gan + \mbox{an }_{\omega, g}(\bar m_1)(\alpha) \bar m_2 \, \forall \alpha \in C(\mathbb{R}^n)
\]
\end{lem}

\begin{pf}
\begin{eqnarray}
&&\lan \bar m_1, \bar \alpha m_2 \gan  \nonumber \\
&=& [m_1, \alpha m_2] + \frac{\omega (\alpha m_2)}{\omega g} [g, m_1] - \frac{\omega m_1}{\omega \alpha} [g, \alpha m_2]  \label{propItemp1}\\
&=& \alpha [m_1, m_2] + m_1(\alpha) m_2 \nonumber \label{propItemp2} \\
&& + \alpha \left( \frac{\omega m_2}{\omega g}[g, m_1] - 
  \frac{\omega m_1}{\omega g}[g, m_2]  \right) \nonumber \label{propItemp3} \\
&& - \frac{\omega m_1 }{\omega g} g(\alpha) m_2 \label{propItemp4} \\
&=& \alpha \lan \bar m_1, \bar m_2 \gan + \left(m_1 - \frac{\omega m_1}{\omega g} m_1 \right) (\alpha) m_2 \nonumber \\
&=& \alpha \lan \bar m_1, \bar m_2 \gan + \mbox{an }_{\omega, g}(m_1)(\alpha) m_2 \nonumber 
\end{eqnarray}
The transition from (\ref{propItemp1}) to (\ref{propItemp4}) uses the same identity applied twice,
$[m_1,\alpha m_2] = \alpha [m_1, m_2] + m_1(\alpha) m_2$ and $[g, \alpha m_2] = \alpha [g, m_2] + g(\alpha) m_2$. The remaining steps are appropriate groupings of terms.
\end{pf}

\begin{lem}
Property II holds:
\[
\mbox{an }_{\omega, g}( \lan \bar m_1, \bar m_2 \gan) = [ \mbox{an }_{\omega, g}(\bar m_1), \mbox{an }_{\omega, g}(\bar m_2)]
\]
\end{lem}

\begin{pf}
Define $\alpha_1 := \frac{\omega m_1}{\omega g}$ and $\alpha_2 = \frac{\omega m_2}{\omega g}$ so that
\begin{eqnarray}
&&\mbox{an }_{\omega,g}\left( \lan \bar m_1, \bar m_2 \gan \right) \nonumber \\
 &=& \mbox{an }_{\omega,g} \left( [m_1,m_2] + \frac{\omega m_2}{\omega g} [g, m_1] - \frac{\omega m_1}{\omega g} m_2 \right) \nonumber \\
&=& [m_1, m_2] + \alpha_2 [g,m_1] - \alpha_1 [g, m_2] \nonumber \\
&& - \frac{1}{\omega g}  \omega \left( [m_1,m_2] + \alpha_2 [g,m_1] - \alpha_1[g,m_2]\right) g \nonumber
\end{eqnarray}
It also holds, for arbitrary vector fields $f_1$, $f_2$ $\in \Gamma T\mathbb{R}^n$, that 
\[
\omega([f_1,f_2]) = f_1(\omega f_2) - f_2(\omega f_1)
\]
so that
\begin{eqnarray} && \omega ([m_1, m_2] + \alpha_2 [g,m_1] - \alpha_1 [g,m_2]) \nonumber \\
&=& m_1(\omega m_2) - m_2 (\omega m_1) + \alpha_2 g(\omega m_1) - \alpha_2 m_1(\omega g) \nonumber \\
&& - \alpha_1 g(\omega m_2) + \alpha_1 m_2(\omega g) \label{develII1}.
\end{eqnarray}
Next, since $m_1\left(\frac{\alpha}{\beta}\right)  = \frac{\beta m_1(\alpha) - \alpha m_1(\beta)}{\beta^2}$ for $\alpha, \beta \in C(\mathbb{R}^n)$, one has
\begin{eqnarray}
m_1(\omega m_2) - \alpha_2 m_1(\omega g) &=& m_1(\omega m_2) - \frac{\omega m_2}{\omega g} m_1(\omega g) \nonumber \\
&=& \omega g \frac{(\omega g) m_1(\omega m_2)- (\omega m_2) m_1(\omega g)}{(\omega g)^2} \nonumber \\
&=& (\omega g) m_1(\alpha_2) \label{develII2}
\end{eqnarray}
Similarly,
\begin{equation}
\label{develII3}
m_2(\omega m_1) - \alpha_1 m_2(\omega g) = (\omega g) m_2(\alpha_1)
\end{equation}
Another expansion gives
\begin{eqnarray}
&&\alpha_2 g(\omega m_1) - \alpha_1 g(\omega m_2) = \alpha_2 g(\omega m_1) - \alpha_2 \alpha_1 g(\omega g) \nonumber \\
&& - \alpha_1 g(\omega m_2) + \alpha_2 \alpha_1 g(\omega g) \nonumber \\
&=& \alpha_2 (g(\omega m_1) - \alpha_1 g(\omega g)) - \alpha_1 (g(\omega m_2) - \alpha_2 g(\omega g)) \nonumber \\
\label{develII4}
&=& \alpha_2(\omega g)g(\alpha_1) - \alpha_1 (\omega g) g(\alpha_2) 
\end{eqnarray}
so that substituting (\ref{develII2}), (\ref{develII3}) and (\ref{develII4}) into (\ref{develII1}) modifies the left-hand side of the identity to be proved in the following way:
\begin{eqnarray}
&&\mbox{an }_{\omega,g}\left(\lan \bar m_1, \bar m_2 \gan \right) = [m_1,m_2] + \alpha_2 [g,m_1] - \alpha_1 [g, m_2] \nonumber \\
&&- \frac{1}{\omega g} ( (\omega g) m_1(\alpha_2) - (\omega g) m_2(\alpha_1)\nonumber \\
&& + \alpha_2 (\omega g) g(\alpha_1) - \alpha_1(\omega g) g(\alpha_2) ) g \nonumber \\
&=& [m_1,m_2] + \alpha_2 [g,m_1] - \alpha_1[g,m_2] \nonumber \\
&& - \left(m_1(\alpha_2) - m_2(\alpha_1) + \alpha_2 g(\alpha_1) - \alpha_1 g(\alpha_2)  \right)g 
\label{develII5}
\end{eqnarray}
Now consider the right-hand side of the identity, namely
\begin{eqnarray}
&&[\mbox{an }_{\omega, g}(\bar m_1), \mbox{an }_{\omega,g}(\bar m_2)] = [m_1  - \alpha_1 g, m_2 - \alpha_2 g] \nonumber \\
&=& [m_1,m_2] - [m_1, \alpha_2 g] - [\alpha_1 g, m_2] - [\alpha_1 g, m_2] + [\alpha_1 g, \alpha_2 g] \nonumber \\
&=& [m_1,m_2] - \alpha_2 [m_1,g] - m_1(\alpha_2) g - \alpha_1 [g, m_2] + m_2(\alpha_1) g \nonumber \\
&& + \alpha_1 \alpha_2 [g,g] + \alpha_1 g(\alpha_2) g - \alpha_2 g(\alpha_1) g \nonumber \\
&=& [m_1,m_2] + \alpha_2 [g,m_1] - \alpha_1 [g,m_2] \nonumber \\
&& + (- m_1(\alpha_2) + m_2(\alpha_1) + \alpha_1 g(\alpha_2) - \alpha_2 g(\alpha_1)) g \label{develII6}
\end{eqnarray}
Comparing (\ref{develII5}) with (\ref{develII6}) shows that
\[
\mbox{an }_{\omega,g}\left( \lan \bar m_1, \bar m_2 \gan\right) = [ \mbox{an }_{\omega,g}(\bar m_1),
\mbox{an }_{\omega,g}(\bar m_2)]
\]
which proves the assertion.
\end{pf}

\subsubsection{Proof of the Jacobi identity}
\begin{lem}
The following identity
\[
\sum_{\mbox{cyclic }i,j,k} \lan m_i, \lan m_j, m_k \gan \gan = 0
\]
holds.
\end{lem}

\begin{pf}
For notation convenience, the following quantities are defined:
\begin{eqnarray*}
\alpha_1:= \frac{\omega m_1}{\omega g} \hspace{1.5cm}\alpha_2:= \frac{\omega m_2}{\omega g}  \hspace{1.5cm}\alpha_3:= \frac{\omega m_3}{\omega g}.
\end{eqnarray*}
Considering the first term of the Jacobi identity and the identity (\ref{equivClass})
\begin{eqnarray}
&&\lan \lan  m_1,  m_2 \gan, m_3 \gan = \nonumber \\
&&[\mbox{an }_{\omg}(\lan  m_1,  m_2 \gan),\mbox{an }_{\omg} (m_3)] \nonumber\\ &+&\left({(\mbox{an}_{\omg}(\lan m_1,m_2 \gan) )} \alpha_3\right)g \nonumber\\&-& \left({\mbox{an }_{\omg}(m_3) }\frac{\omega \lan m_1,m_2 \gan}{\omega g}\right) g. \label{jacob1}
\end{eqnarray}
By using (\ref{equivClass}) for $\lan m_1,m_2 \gan$, we get
\begin{eqnarray}
&&\omega \lan  m_1, m_2  \gan =\omega ( [\mbox{an }_{\omg}(m_1),\mbox{an }_{\omg}(m_2)]  \nonumber \\
&&+({\mbox{an }_{\omg}(m_1)}\alpha_2 - {\mbox{an }_{\omg}(m_2)}\alpha_1)g)\nonumber\\
&=& \omega[\mbox{an}_{\omg}(m_1),\mbox{an }_{\omg}(m_2)] + ({\mbox{an }_{\omg}(m_1)}\alpha_2 \nonumber \\
&& - {\mbox{an }_{\omg}(m_2)}\alpha_1) \omega g\nonumber\\
&=&0 + ({\mbox{an }_{\omg}(m_1)}\alpha_2 - {\mbox{an }_{\omg}(m_2)}\alpha_1) \omega g.\label{jacob_int}
\end{eqnarray}
Substituting (\ref{jacob_int}) in (\ref{jacob1}) gives with $i=1, j=2, k=3$
\begin{eqnarray*}
&&\lan \lan m_i,m_j \gan,m_k\gan \nonumber \\
&=&[[\mbox{an }_{\omg}(m_i),\mbox{an }_{\omg}(m_j)] \mbox{an }_{\omg}(m_k)] \\&&+\left({[\mbox{an }_{\omg}(m_i), \mbox{an }_{\omg}(m_j)]}(\alpha_k)\right.\\&&\left.-{\mbox{an }_{\omg}(m_k) } \frac{({\mbox{an }_{\omg}(m_i)} \alpha_j - {\mbox{an }_{\omg}(m_j)}(\alpha_i)) \omega g} {\omega g}\right)g
\end{eqnarray*}
\begin{eqnarray*}
&=&[[\mbox{an }_{\omg}(m_i),\mbox{an }_{\omg}(m_j)] \mbox{an }_{\omg}(m_k)] \\&&+({\mbox{an }_{\omg}(m_i} {\mbox{an }_{\omg}(m_j)}\alpha_k-{\mbox{an }_{\omg}(m_j)} {\mbox{an }_{\omg}(m_i)}\alpha_k\\&&-{\mbox{an }_{\omg}(m_k) }{\mbox{an }_{\omg}(m_i)}\alpha_j + {\mbox{an }_{\omg}(m_k) }{\mbox{an }_{\omg}(m_j)}\alpha_i)g.
\end{eqnarray*}
It is then straightforward to notice that 
a circular summation of the previous expression over the indices $i,j,k$ yields zero, that is,
\begin{eqnarray*}
\sum_{\mbox{cyclic }i,j,k} \lan \lan m_i,m_j \gan,m_k\gan = 0 
\end{eqnarray*}
which is the Jacobi identity.
\end{pf}

\subsection{Lie Algebroid Isomorphism}

\begin{prop}
The algebroids of Sections \ref{Algebroid1} and \ref{Algebroid2} are isomorphic in the sense
that there exists a one-to-one correspondance between $\mathcal O$ - projectable vectorfields
and corresponding line bundle in the $g, \omega$-quotient bundle.
\end{prop}

\begin{pf}The right-hand-side of (\ref{RHS2}) is the same as the right-hand-side
of Property (II) of the algebroid of the groupoid. Therefore, if one gives
two $\mathcal O$ - projectable vectorfields ${\tilde m}_1$ and ${\tilde m}_2$, then one simply defines corresponding line bundles as
$\{ {\tilde m}_1 + \alpha g, \forall \alpha : \mathbb{R}^N \rightarrow \mathbb{R} \}$ and 
$\{ {\tilde m}_2 + \alpha g, \forall \alpha : \mathbb{R}^N \rightarrow \mathbb{R} \}$  
 for which ${\tilde m}_1$ and ${\tilde m}_2$ are used as representatives. Then $\pi( \langle {\tilde m}_1, {\tilde m}_2 \rangle) = \pi([{\tilde m}_1, {\tilde m}_2]) = [\pi{\tilde m}_1, \pi{\tilde m}_2] = [{\bar m}_1, {\bar m}_2]$.
 Reciprocally, suppose that two line bundles are given a priori, namely $\{m_1 + \alpha g, \forall \alpha : \mathbb{R}^N \rightarrow \mathbb{R} \}$, and 
$\{ m_2 + \alpha g, \forall \alpha : \mathbb{R}^N \rightarrow \mathbb{R} \}$ and
compute $\bar m_1 = \pi(m_1) = \mbox{Pr}(\Phi_*(m_1)) $ and $\bar m_2 = \pi(m_2) = \mbox{Pr}(\Phi_*(m_2))$ so that after setting
$${\tilde m}_1 = (\Phi_*)^{-1} \begin{pmatrix} {\bar m}_1  \cr 0 \end{pmatrix}  \qquad {\tilde m}_2 = (\Phi_*)^{-1}\begin{pmatrix} {\bar m}_2 \cr 0 \end{pmatrix}.$$ one notices that because of the zero inserted in the last component, the vectorfields ${\tilde m}_1$ and ${\tilde m}_2$ are  $\mathcal O$ - projectable and therefore satisfy $\pi([{\tilde m}_1, {\tilde m}_2]) = [\bar m_1, \bar m_2] = \pi(\langle {\tilde m}_1, {\tilde m}_2 \rangle)$
Because by construction of ${\bar m}_1$ and ${\bar m_2}$, it is true that  $[{\bar m}_1, {\bar m}_2] = \pi(\langle m_1, m_2 \rangle)$, this also means that ${\tilde m}_1$ belongs to the line bundle generated by $m_1$, and ${\tilde m}_2$ belongs to the line bundle generated by $m_2$. The arbitrariness of $m_1$ and $m_2$ within their respective line bundles shows that the construction of $\tilde m_1$ and $\tilde m_2$ does not depend on the representatives $m_1$ and $m_2$ chosen. 

Therefore, a one-to-one correspondance between $\mathcal O$ - projectable vectorfields and corresponding line bundles is established. The elements of one set (the $\mathcal O$ - projectable vectorfields $\tilde m_1$ or $\tilde m_2$) or the other (the line bundles $\{m_1 + \alpha g, \forall \alpha: \mathbb{R}^n \rightarrow \mathbb{R}\}$ or
$\{m_2 + \alpha g, \forall \alpha: \mathbb{R}^n \rightarrow \mathbb{R}\}$) are distinguished by the vectorfields $\bar m_1$ and $\bar m_2$ to which they map in $T \mathcal O$. \end{pf}

\section{Application to Feedback Linearization}

\label{AppFbkLin}

\subsection{Algorithm using Algebroid I}

This algorithm is described in \cite{Mullhaupt} and is summarized hereafter.
It consists of two phases. The first phase reduces the number of coordinates using diffeomorphisms of the Lie groupoid,  keeping track of their inverses. The linearizing output is computed using the chain of inverses of the target maps during the second phase.

\subsubsection{Phase 1}
\begin{itemize}
\item \emph{Initialisation}: $f_0 := f$, $g_0 := g$ and define $\mbox{an}_{\pi,0}$ using a diffeomorphism $\Phi_0$ such that $\mbox{an}_{\pi,0}(g_0) = 0$.
\item \emph{Induction}: \begin{eqnarray}
 f_{i+1} &=& \mbox{an}_\pi(f_i)  \nonumber \\
g_{i+1} &=& \mbox{an}_{\pi,i}([f_i,g_i]) \nonumber 
\end{eqnarray}
and choose $\omega_{i+1}$ such that it is integrable (or exact) such that $\omega_{i+1} g_{i+1} \neq 0$ and construct a diffeomorphism $\Phi$ associated with the groupoid and defining $\mbox{an}_{\pi,i+1}$ such that $\mbox{an}_{\pi,i+1}(g_{i+1}) = 0$.
\item \emph{Termination}: Stop when $i=n-1$.
\end{itemize}

\subsubsection{Phase 2}
The linearizing output is obtained using the chain of inverses of the target maps
\[
z = \Phi_0^{-1} \circ \Phi_1^{-1} \circ \ldots \circ \Phi_{n-1}^{-1} (x_1) 
\]
where $x_1$ stands for the unique state of the last iteration.

\subsection{Algorithm using Algebroid II}

\subsubsection{Phase 1}
This algorithm is described in \cite{Willson} without the formalism of Lie algebroids and groupoids.
\begin{itemize}
\item \emph{Initialisation}: $f_0 := f$, $g_0 := g$ and choose $\omega_0$ integrable (or exact) such that $\omega_0 g_0 \neq 0$.
\item \emph{Induction}: 
\begin{eqnarray}
f_{i+1} &:=& \mbox{an}_{\omega_i,g_i}(f_i) \nonumber \\
 g_{i+1} &:=& \mbox{an}_{\omega_i}([f_i,g_i]) \nonumber 
\end{eqnarray}
Choose $\omega_{i+1}$ integrable (or exact) such that $\omega_{i+1} g_{i+1} \neq 0$.
\item \emph{Termination}: Stop when $i = n-1$.
\end{itemize}

\subsubsection{Phase 2}
\label{Phase2AlgebroidII}
The second phase constructs the linearizing output using the $1$-forms $\omega_i$  used in the first phase:
\begin{itemize}
\item \emph{Initialisation}: $\nu_{n-1} := \omega_{n-1}$
\item \emph{Induction}: \[ \nu_{n-(i+1)} := \nu_{n-i} - \frac{\nu_{n-i} \, g_{n-(i+1)}}{\omega_{n-(i+1)} g_{n-(i+1)}} \omega_{n-(i+1)} \]
\item \emph{Termination}: Stop when $i = n-1$.
\end{itemize}

\section{Polynomial Automorphisms and the Jacobian Conjecture}

Key to all algorithms and properties of the previous sections is the construction of the 1-forms $\omega_i$. The choice of exact forms for which $\omega_i g_i$ are constants  and those that cancel $g_i$ play a fundamental role in the construction of the inverse of a polynomial automorphism as it will be shown in this section through an example.

\label{AppPolAuto}

\subsection{Example}

The polynomial vectorfield $f$ is given by its components $f = \left(\begin{array}{ccc} f_1 & f_2 &f_3 \end{array} \right)^T$ as
\begin{eqnarray}
f_1 &=&
 \frac{x_3^4}{2}+x_2
   x_3^2+\frac{x_3^2}{2} 
 +\frac{x_3}{2}+\frac{x_1^2}{2}
+
   \frac{x_2^2}{2}+\frac{x_1}{2}+\frac{x_2}{2} \nonumber \\
f_2 &=&-4 x_1 x_3^7+2 x_3^7-4 x_1 x_3^5-12 x_1 x_2 x_3^5+6 x_2
   x_3^5 \nonumber \\
&& -2 x_3^5-5 x_1 x_3^4+\frac{5 x_3^4}{2}-4 x_1^3
   x_3^3-6 x_1^2 x_3^3-12 x_1 x_2^2 x_3^3 \nonumber \\
&&+6 x_2^2 \nonumber \\
f_3 &=& x_3^3-2 x_1 x_3^3-8 x_1 x_2 x_3^3-4 x_2 x_3^3-x_1
   x_3^2 \nonumber \\
&& -6 x_1 x_2 x_3^2+3 x_2 x_3^2 \nonumber \\
&& -\frac{x_3^2}{2}-4
   x_1 x_2^3 x_3+2 x_2^3 x_3-4 x_1 x_2^2 x_3-2 x_2^2
   x_3 \nonumber \\
&& -x_1 x_3-4 x_1^3 x_2 x_3-6 x_1^2 x_2 x_3-2 x_1
   x_2 x_3+\frac{x_3}{2} \nonumber \\
&& -x_1^3-\frac{3 x_1^2}{2}-x_1
   x_2^2+\frac{x_2^2}{2}-\frac{x_1}{2}-x_1
   x_2-\frac{x_2}{2} \nonumber
\end{eqnarray}
and the $g$ vectorfield is 
\begin{eqnarray}
g &=&  \left(
\begin{array}{ccc}
 0 & -2 x_3 & 1 \\
\end{array}
\right)^T \nonumber 
\end{eqnarray}

The polynomial vectorfields $f$ and $g$ can be understood as polynomial derivations $f= \sum_i f_i \frac{\partial}{\partial x_i}$ and $g = \sum_i g_i \frac{\partial}{\partial x_i}$ \cite{vanEssen}.

\subsection{Algorithm with Algebroid II}

\subsubsection{Phase 1}

The indices of $f$ now relate to the iteration number of the algorithm (and not to its components). Hence set $f_0 = f$ and $g_0 = g$.
The 1-form
\[
\omega_0 = (2 x_3^2+2 x_2) dx_2 + (x_3^3+4 x_2 x_3+1) dx_3 
\]
is such that  $ \omega_0 \, g_0 = 1$ and is exact since $\omega_0 = d(x_3^4+2 x_2 x_3^2+x_3+x_2^2)$.
This will be used to define the first anchor
\[
\mbox{an }_{\omega_0, g_0}(m) = m - \frac{\omega_0 \, m}{\omega_0 \, g_0} g_0
\]
A direct computation gives
\begin{eqnarray}
g_1 &=& \mbox{an }_{\omega_0, g_0} ([f_0, g_0]) \nonumber \\
&=&
\left(
\begin{array}{c}
 -\frac{1}{2} \\
 4 x_1 x_3^3-2 x_3^3+4 x_1 x_2 x_3-2 x_2 x_3+x_1-\frac{1}{2} \\
 -2 x_1 x_3^2+x_3^2-2 x_1 x_2+x_2 \\
\end{array}
\right) \nonumber 
\end{eqnarray}
and $f_1 = f_0$.
Selecting the trivial exact 1-form
\[
\omega_1 = dx_1
\]
leads to the second iteration which is 
\begin{eqnarray}
g_2 &=& \mbox{an }_{\omega_2, g_2} ([f_1, g_1]) =  \nonumber \\
&=& \left(
\begin{array}{ccc}
 0 & -4 x_3^3-4 x_2 x_3-1 & 2 \left(x_3^2+x_2\right) \nonumber 
\end{array}
\right)^T  \nonumber
\end{eqnarray}
Choose $\omega_2 = d x_2$ so that 
\[
\omega_2 \, g_2 = -4 x_3^3-4 x_2 x_3-1
\]
this will be the integrating factor of the $1$-form $\nu_0$ constructed in Phase 2.

\subsubsection{Phase 2}

Applying the iteration scheme of Section \ref{Phase2AlgebroidII} gives
\begin{eqnarray}
\nu_2 &=& \omega_2 = d x_2 \nonumber \\
\nu_1 &=&  (8 x_1 x_3^3-4 x_3^3+8 x_1 x_2 x_3-4 x_2 x_3+2
   x_1-1)dx_1 + dx_2 \nonumber \\
\nonumber \\
\nu_0 &=& (8 x_1 x_3^3-4 x_3^3+8 x_1 x_2 x_3-4 x_2 x_3+2 x_1-1) dx_1 \nonumber \\
&& + (4
   x_3^3+4 x_2 x_3+1) dx_2 + (8 x_3^4+8 x_2 x_3^2+2 x_3) dx_3 \nonumber 
\end{eqnarray}

Integrating the exact form $\frac{1}{\omega_2 \, g_2} \nu_0$ leads to the linearizing output
\begin{eqnarray}
y = \int \frac{1}{\omega_2 \, g_2} \nu_0 =  x_1 - x_1^2 - x_2  - x_3^2 \nonumber
\end{eqnarray}

\subsection{Algorithm with Algebroid I}

\subsubsection{Phase 1}

Set $f_0=f$ and $g_0=g$.
The polynomial morphism
\[
\Phi_0 : x \rightarrow \left(
\begin{array}{c}
 x_1 \\
 x_3^2+x_2 \\
 x_3^4+2 x_2 x_3^2+x_3+x_2^2 \\
\end{array}
\right)
\]
 admits the inverse
\[
\Phi_0^{-1} : z \rightarrow
\left(
\begin{array}{c}
 z_1 \\
 -z_2^4+2 z_3 z_2^2+z_2-z_3^2 \\
 z_3-z_2^2 \\
\end{array}
\right)
\]
so that the anchor 
\[
\mbox{an }_{\pi,0}(m) =
 \mbox{Pr } \Phi_{*,0}(m)
\] is defined such that $\mbox{an }_{\pi,0}(g_0) = 0$.
Then 
\begin{eqnarray}
f_1 &=& \mbox{an }_{\pi, 0} (f_0)  \nonumber  \\
&=&
\left(
\begin{array}{c}
 \frac{z_1^2}{2}+\frac{z_1}{2}+\frac{z_2}{2}+\frac{z_3}{2} \\
 -z_1^3-\frac{3 z_1^2}{2}-z_2 z_1-z_3
   z_1-\frac{z_1}{2}-\frac{z_2}{2}+\frac{z_3}{2} \\
\end{array}
\right)
 \nonumber \\
g_1 &=& \mbox{an }_{\pi, 0} ([f_0,g_0]) = 
\left(
\begin{array}{c}
 -\frac{1}{2} \\
 z_1-\frac{1}{2} \\
\end{array}
\right)\nonumber 
\end{eqnarray}
Select the second polynomial morphism as
\begin{eqnarray}
\Phi_1: z \rightarrow \left( \begin{array}{c}
z_1 + z_1^2 + z_2 \nonumber \\
z_1 - z_1^2 - z_2 \nonumber \end{array} \right) \nonumber
\end{eqnarray}
with polynomial inverse
\begin{eqnarray} \kern-0.5cm
\Phi_1^{-1}: w \rightarrow 
\left(
\begin{array}{c}
 \frac{1}{2} \left(w_1+w_2\right) \\
 \frac{1}{4} \left(-w_1^2-2 w_2 w_1+2 w_1-w_2^2-2 w_2\right) \\
\end{array}
\right) \label{PhiInv1} 
\end{eqnarray}
defining the second anchor
\[
\mbox{an }_{\pi,1}(m) = \mbox{Pr } \Phi_{*,1}(m)
\]
with the property that $\mbox{an }_{\pi,1} (g_1) = 0$.
The linearizing output is $w_1$.

\subsubsection{Phase 2}

Phase 2 consists in expressing $w_1$ through the successive polynomial-inverse maps:

\begin{eqnarray}
y &=& \Phi_0^{-1} ( \Phi_1^{-1} (w_1)) 
= \Phi_0^{-1} (z_1 - z_2^2 - z_2) \nonumber \\
&=& x_1 - x_1^2 - x_2 - x_3^2 \nonumber
\end{eqnarray}

\subsection{Relation to the Jacobian Conjecture}

Setting 
\begin{eqnarray}
\Phi : x  \rightarrow  \left(
\begin{array}{c}
y \\
L_f y \\
L_f^2 y
\end{array}
\right) = \left(\begin{array}{c}
 x_1 - x_1^2 - x_2 - x_3^2 \\
 x_1^2+x_1+x_3^2+x_2  \\
 x_3^4+2 x_2 x_3^2+x_3+x_2^2 
\end{array} \right) \label{PhiMap}
\end{eqnarray}
gives a polynomial morphism $\Phi : x \rightarrow \Phi(x)$. Extending the $\Phi_1$ map obtained in Phase 2 with $z_3 \rightarrow z_3$ and changing notations using $x$ instead of $z$ gives the polynomial morphism
\begin{eqnarray}
\Psi : x \rightarrow \left( \begin{array}{c} x_1 + x_1^2 +x_2 \\
x_1 - x_1^2 - x_2 \\
x_3
\end{array} \right) \nonumber 
\end{eqnarray}
with inverse given as (\ref{PhiInv1}) with $w$ replaced by $x$ and with last component $x_3$.
It is then straightforward to show that $\Psi^{-1} \circ \Phi_0^{-1}$ is the inverse map of $\Phi$ defined in (\ref{PhiMap}). 

Associated with any polynomial automorphism, one can construct a dynamical system $\dot x = f(x) + g(x) u$ which is feedback linearizable using the polynomial automorphism. With $n=3$ this would be $\dot z_1=z_2$, $\dot z_2 = z_2$, $\dot z_3 = u$,
and determine the associated $f$ and $g$ using the polynomial morphism. Then proceed as described with $f$ and $g$ given above. The example was constructed using a particular class of tame polynomial automorphisms. 

\section{Conclusion}

The algebroids given in Section \ref{Algebroid1} and \ref{Algebroid2} have different anchors and can be used to give two iteratives schemes to compute the linearizing output of nonlinear affine in the input single-input system.  The algebroids were shown to  satisfy the Jacobi identity and all properties required.  Key in establishing this result is the fact that $\omega$ appearing in (\ref{equivClass}) is an integrable 1-form. Using the two algebroids an example using polynomial automorphisms instead of diffeomorphisms illustrated the theory. 
The convergence and computation of the inverse polynomial map
hinged on the construction of exact forms in the intermediate steps of the algorithm. An algorithm for a class of tame polynomial automorphisms was used for generating the example and will be described elsewhere. 


\end{document}